\documentclass[11pt,bezier]{article}
\usepackage{amsmath,amssymb,amsfonts,euscript}
\newpage
\newcommand{\be}{\begin{equation}}
\newcommand{\ee}{\end{equation}}
\newtheorem{theorem}{Theorem}[section]

\newtheorem{example}{Example}[section]
\newtheorem{definition}{Definition}[section]
\newtheorem{remark}{Remark}[section]
\newtheorem{proposition}{Proposition}[section]
\renewcommand{\theequation}{\arabic{section}.\arabic{subsection}.\arabic{equation}}
\title{\bf\Large A new structure for analyzing discrete scale\\ invariant processes:  Covariance and Spectra}
\vspace{-1cm}
\author
{ N. Modarresi\,\,\,\, and \,\,\,  S. Rezakhah\thanks{ Faculty of
Mathematics and Computer Science, Amirkabir University of
Technology, 424 Hafez Avenue, Tehran 15914, Iran. E-mail:
namomath@aut.ac.ir (Modarresi), rezakhah@aut.ac.ir (Rezakhah).} }
\date{}
\begin{document}
\maketitle

\begin{abstract}

Improving the efficiency of discrete time  scale invariant (DSI)
processes, we consider some flexible sampling of a continuous time
DSI process $\{X(t), t\in{\bf R^+}\}$ with scale $l>1$, which is
 in correspondence to some  multi-dimensional self-similar process.  So we  consider
 $q$ samples at arbitrary points ${\bf s}_0,
{\bf s}_1, \ldots, {\bf s}_{q-1}$ in interval $[1, l)$ and proceed
in the intervals $[l^n, l^{n+1})$ at points $l^n{\bf s}_0,l^n{\bf
s}_1, \ldots, l^n{\bf s}_{q-1}$, $n\in {\bf Z}$. So we study
 an embedded DT-SI
process  $W(nq+k)=X(l^n{\bf s}_k)$, $q\in {\bf N}$, $k= 0, \ldots,
q-1$, and its multi-dimensional self-similar counter part  $V(n)=\big(V^0(n), \ldots, V^{q-1}(n)\big)$
where $V^k(n)=W(nq+k)$. We study spectral representation of such process and obtain its spectral
density matrix.
Finally by imposing wide sense Markov property on $W(\cdot)$ and $V(\cdot)$, we show that the spectral density matrix of $V(\cdot)$
 can be characterized by $\{R_j(1), R_j(0), j=0, \ldots, q-1\}$ where $R_j(k)=E[W(j+k)W(j)]$.\\
%

{\it AMS 2010 Subject Classification:} 60G18, 62M15, 60J05.\\

{\it Keywords:} Discrete scale invariance; Wide sense Markov; Multi-dimensional self-similar process; Spectral density matrix; Brownian motion.
\end{abstract}

\section{Introduction}
The concept of stationarity and self-similarity are used as a
fundamental property to handle many natural phenomena. Lamperti
transformation defines a one to one correspondence between
stationary and self-similar processes. Discrete scale invariant
(DSI) processes can be defined as
the Lamperti transform of periodically correlated ones. Many
critical systems, like statistical physics, textures in geophysics,
network traffic and image processing can be interpreted by these
processes \cite {b1}. Fourier transform is known as a suited
representation for stationary processes.  Using Mellin transform, a
harmonic like representation of self-similar process is introduced \cite {f1}.
Self-similar Markov process which has Markov property and self-similarity
are involved in various parts of probability theory, such as branching processes
and fragmentation theory \cite {c1}. Gladyshev in \cite{g1} introduced the spectral
representation of correlation matrix of multi-dimensional stationary random sequences
and found a relation between them and periodically correlated (PC) processes.

In our previous work \cite{m2} we studied spectral analysis of a sequence of observations which are sampled at some special points, $\alpha^k$, $k\in {\Bbb Z}$ of a DSI process with some scale $l=\alpha^T$, $T\in {\Bbb N}$. So that one could study such processes in spectral domain.
By such sampling, and by imposing wide sense Markov property, we
provided a discrete time scale invariant sequence which is Markov in the wide sense and
found a closed formula for its covariance function and spectral density matrix of corresponding multi-dimensional self-similar process. The above sampling scheme had much restrictions which could dismiss lots of information between sample points $\alpha^k$.
In this paper we consider some flexible sampling scheme which enables one to have samples at arbitrary points ${\bf s}_0, {\bf s}_1, \ldots, {\bf s}_{q-1}$ in the first scale interval and follow sampling at corresponding points $\{l^n{\bf s}_i, n\in\Bbb N, i=0, \ldots, q-1\}$ in the other scale intervals.
  This sampling scheme  provide
  a corresponding  multi-dimensional self-similar process as a platform to extend analytic property of discrete time periodically correlated Markov process to such DSI sequences.
For such study we need to consider a new approach for modeling this sequence via defining some embedded process. So this method enables one to study the behavior of the process specially in spectral domain and have a better description of the underlying DSI process at all arbitrary points.

Let $X(\cdot)$ be a DSI process with scale $l>1$. By sampling at arbitrary points in the scale
intervals $[l^{n-1}, l^{n})$, $n\in {\Bbb N}$, we provide $\{X(l^{n-1}{\bf s}_k), n\in {\Bbb N};
1\leq {\bf s}_0<\cdots <{\bf s}_{q-1}<l\}$ as a sequence of DSI.
Then we define an embedded scale invariant process as $\{W(nq+k)\equiv X(l^n{\bf s}_k); n\in{\Bbb W}, k=0, \ldots, q-1\}$ where $\Bbb W=\{0, 1, \ldots\}$ and its corresponding multi-dimensional embedded
self-similar process as $\{V(n), n\in {\Bbb W} \}$, where
$V(n)=\big(V^0(n), \ldots, V^{q-1}(n)\big)$, and $\{V^i(n)\equiv W(nq+i), n\in{\Bbb W}, i=0, \ldots, q-1\}$,
for fixed $q\in {\Bbb N}$.
We investigate covariance structure and spectral density matrix of $\{U(l^n), n\in {\Bbb W} \}$, where
$U(l^n)=\big(U^0(l^n), \ldots, U^{q-1}(l^n)\big)$, and  $\{U^i(l^n)\equiv X(l^n{\bf s}_i), n\in{\Bbb W}, i=0, \ldots, q-1\}$, for fixed $q\in {\Bbb N}$, when it is Markov in the
wide sense and is called multi-dimensional self-similar Markov process.

This paper is organized as follows. In section 2, we present multi-dimensional stationary, PC,
self-similar and DSI processes and define them in discrete time. We also review the properties of
Lamperti transformation in this section.
Section 3 is devoted to the structure of the multi-dimensional self-similar process resulting from the above method of
sampling. We define embedded scale invariant process and corresponding
multi-dimensional embedded self-similar process and characterize the
spectral density matrix of it in this section. We also characterize
covariance function and spectral density matrix of the multi-dimensional
self-similar and embedded self-similar Markov processes in section 4 which is
characterized by $\{R_j(1), R_j(0), j=0, \ldots, q-1\}$ where
$R_j(k)=E[W(j+k)W(j)]$. We simulate
and analyze Simple Brownian Motion (SBM) as a DSI Markov and its
corresponding multi-dimensional processes. Throughout this paper we study above
mentioned processes in the discrete time with some clarified parameter spaces.
So we omit the term discrete time in the rest of this paper.

\renewcommand{\theequation}{\arabic{section}.\arabic{equation}}
\section{Theoretical framework}
\setcounter{equation}{0}
In this section we review the structure of covariance function and spectral density matrix of multi-dimensional stationary processes. The definitions of self-similar, scale invariant processes in discrete time and in the wide sense are presented. Also Lamperti transformation is defined and its properties are studied.

\subsection{Stationary and multi-dimensional stationary processes}
\begin{definition}
A process $\{Y(t),t\in {\Bbb R}\}$ is said to be stationary, if for any $\tau\in {\Bbb R}$
\be \{Y(t+\tau), t\in {\Bbb R} \}\stackrel{d}{=}\{Y(t),t\in {\Bbb R}\}\ee
where $\stackrel{d}{=}$ is the equality of all finite-dimensional distributions.

If $(2.1)$ holds for some $\tau\in {\Bbb R}$, the process is said to be periodically correlated. The smallest of such $\tau$ is called period of the process.
\end{definition}

By Rozanov \cite {r1}, if $Y(t)=\{Y^k(t)\}_{k=1,\ldots,n}$ be an $n$-dimensional stationary process, then
\be Y(t)=\int e^{i\lambda t}\phi(d\lambda)\ee
is its spectral representation, where $\phi=\{\varphi_k\}_{k=1,\ldots,n}$ and $\varphi_k$ is the random spectral measure associated with the $k$th component $Y^k$ of the $n$-dimensional process $Y$. Let
$$B_{kr}(\tau)=E[Y^k(\tau+t)\overline{Y^r(t)}], \hspace{7mm}k,r=1,\ldots,n$$
and $B(\tau)=[B_{kr}(\tau)]_{k,r=1,\ldots,n}$ be the correlation matrix of $Y$.
The components of the correlation matrix of the process $Y$ can be represented as
\be B_{kr}(\tau)=\int e^{i\lambda \tau}F_{kr}(d\lambda), \hspace{7mm}k,r=1,\ldots,n\ee
where for any Borel set $\Delta$, $F_{kr}(\Delta)=E[\varphi_k(\Delta)\overline{\varphi_r(\Delta)}]$  are the complex valued set functions which are $\sigma$-additive and have bounded variation. For any $k,r=1,\ldots,n$, if the sets $\Delta$ and $\Delta'$ do not intersect, $E[\varphi_k(\Delta)\overline{\varphi_r(\Delta')}]=0$.
For any interval $\Delta=(\lambda_1,\lambda_2)$ when $F_{kr}(\{\lambda_1\})=F_{kr}(\{\lambda_2\})=0$ the following relation holds
\be F_{kr}(\Delta)=\frac{1}{2\pi}\int_{\Delta}\sum_{\tau=-\infty}^{\infty}B_{kr}(\tau)e^{-i\lambda \tau}d\lambda\ee
$$=\frac{1}{2\pi}B_{kr}(0)[\lambda_2-\lambda_1]+\lim_{T\rightarrow\infty}\frac{1}{2\pi}\sum_{0<|\tau|\leqslant T}
B_{kr}(\tau)\frac{e^{-i\lambda_2 \tau}-e^{-i\lambda_1 \tau}}{-i\tau}$$
in the discrete parameter case, and
$$F_{kr}(\Delta)=\lim_{a\rightarrow\infty}\frac{1}{2\pi}\int_{-a}^{a}\frac{e^{-i\lambda_2
\tau}-e^{-i\lambda_1 \tau}}{-i\tau}B_{kr}(\tau)d\tau$$
in the continuous parameter case.

\subsection{Discrete time scale invariant processes}
\begin{definition}
A process $\{X(t),t\in {\Bbb R^+}\}$ is said to be self-similar of
index $H>0$, if for any ${\lambda}>0$
\be \{\lambda^{-H}X(\lambda t), t\in {\Bbb R}\}\stackrel{d}{=}\{X(t), t\in {\Bbb R}\}.\ee

The process is said to be {\em DSI} of index $H$ and scaling factor ${\lambda}_0>0$ or {\em (H,${\lambda}_0$)-DSI}, if $(2.5)$ holds for $\lambda=\lambda_0$.
\end{definition}

As an intuition, self-similarity refers to an invariance with respect to any dilation factor. However, this may be a too strong requirement for capturing in situations that scaling properties are only observed for some preferred dilation factors.

\begin{definition}
A process $\{X(k),k\in {\hat{T}}\}$ is called discrete time self-similar process with parameter space $\hat{T}$,
where $\hat{T}$ is any subset of countable distinct points of positive real numbers, if for any $k_1, k_2 \in \check{T}$
\be \{X(k_2)\}\stackrel{d}{=}(\frac{k_2}{k_1})^H\{X(k_1)\}.\ee
The process $X(\cdot)$ is called discrete time scale invariant with scale $l>0$ and parameter space $\hat{T}$, if for any $k_1, k_2=lk_1 \in \hat{T}$, $(2.6)$ holds.
\end{definition}

\begin{remark}
If the process $\{X(t),t\in {\Bbb R^+}\}$ is {\em DSI} with scale $l>1$. Then by sampling of the process at points of set
$$\check{T}=\{l^{n}{\bf s}_j, n\in {\Bbb Z}, 1\leq {\bf s}_0<\cdots <{\bf s}_{q-1}<l\}$$
for $j=0, \ldots, q-1$, we have $X(\cdot)$ as a scale invariant process with parameter space $\check{T}$ and scale $l$.
If we consider sampling of $X(\cdot)$ at points
$$\tilde{T}=\{l^{n}{\bf s}_j, n\in {\Bbb Z}\ \ \text{for fixed}\ \ 1\leq{\bf s}_j<l\},$$
then $X(\cdot)$ is a self-similar process with parameter space $\tilde{T}$.
\end{remark}

\noindent Based on the definition of wide sense self-similar process
presented in \cite{n1}, we present the following definition.

\begin{definition}
A random process $\{X(k),k\in\hat T\}$ is called self-similar in the wide sense with index $H>0$ and with parameter space $\hat T$, where $\hat T$ is any subset of distinct countable points of positive real numbers, if for all $k, k_1\in\hat T$ and all $c>0$, where $ck, ck_1\in\hat T$\\

$(i)\,\,\ E[X^2(k)]<\infty$,

$(ii)\,\,E[X(ck)]=c^HE[X(k)]$,

$(iii)\,\, E[X(ck)X(ck_1)]=c^{2H}E[X(k)X(k_1)]$.\\ \\
If the above conditions hold for some fixed $c=c_0$, then the process is called scale invariant in the wide sense with scale $c_0$.
\end{definition}

Throughout this paper we are dealt with wide sense self-similar and wide sense scale invariant process, and for simplicity we omit the term "in the wide sense" hereafter.

\subsection{Lamperti transformation}
The Lamperti transformation provides a bijection between self-similar and stationary processes and also scale invariant and periodically correlated processes in discrete time.
\begin{definition}
The Lamperti transform with positive index $H$, denoted by ${\cal L}_H$ operates on a random process
$\{Y(t),t\in {\Bbb R}\}$ as
\be {\cal L}_HY(t)=t^HY(\ln t)\ee
and the corresponding inverse Lamperti transform ${\cal L}^{-1}_H$ on process $\{X(t), t\in {\Bbb R^+}\}$
acts as
\be{\cal L}^{-1}_HX(t)=e^{-tH}X(e^t).\ee
\end{definition}

\begin{remark}
If $\{Y(t),t\in {\Bbb R}\}$ is stationary process, its Lamperti transform  $\{{\cal L}_HY(t), t\in {\Bbb R^+}\}$ is self-similar. Conversely if $\{X(t),t\in {\Bbb R^+}\}$ is self-similar process, its inverse Lamperti transform $\{{\cal L}^{-1}_HX(t),t\in {\Bbb R}\}$ is stationary.
\end{remark}

\begin{remark}
If $\{X(t), t\in {\Bbb R^+}\}$ is {\em $(H,e^T$)-DSI} then ${\cal L}^{-1}_HX(t)=Y(t)$ is periodically correlated with period $T>0$. Conversely if $\{Y(t), t\in {\Bbb R}\}$ is periodically correlated with period $T$ then ${\cal L}_HY(t)=X(t)$ is {\em $(H,e^T$)-DSI}.
\end{remark}

\begin{remark}
If $X(\cdot)$ is a self-similar process with parameter space $\check{T}=\{l^n{\bf s}_i,i=0,1,\ldots, q-1; n\in\Bbb W\}$, then its stationary counterpart $Y(\cdot)$ has parameter space $\tilde{T}=\{n\ln l+\ln{\bf s}_i,i=0,1,\ldots, q-1; n\in\Bbb W\}$
$$X(l^{n}{\bf s}_i)
={\cal L}_HY(l^{n}{\bf s}_i)=(l^{n}{\bf s}_i)^HY(n\ln l+\ln{\bf s}_i).$$
Also it is clear by the following relation that if $X(\cdot)$ is a scale invariant process with scale $l$ and parameter space $\check{T}=\{l^{n}{\bf s}_i,i=0,1,\ldots, q-1; n\in\Bbb W\}$, then $Y(\cdot)$
is a discrete time periodically correlated process with period $q$ and parameter space $\tilde{T}=\{n\ln l+\ln{\bf s}_i,
i=0,1,\ldots, q-1; n\in\Bbb W\}$
$$Y(n\ln l+\ln{\bf s}_i)={\cal L}^{-1}_{H}X(n\ln l+\ln{\bf s}_i)=(l^{n}{\bf s}_i)^{-H}X(l^{n}{\bf s}_i).$$
\end{remark}

\renewcommand{\theequation}{\arabic{section}.\arabic{equation}}
\section{Structure of the process}
\setcounter{equation}{0}
In this section we introduce a new method for flexible sampling of a DSI process with scale $l>1$,
which provide sampling at arbitrary points in the interval $[1, l)$ and at multiple $l^n$
of such points in the intervals $[l^n, l^{n+1})$, $n\in {\Bbb N}$. Based on such a sequence of DSI process we define an embedded scale invariant process which follow this sequence. Then by introducing the corresponding multi-dimensional embedded self-similar process in the wide sense we provide a platform to characterize spectral representation and spectral density of the process.
Finally in Theorem 3.1 we find harmonic like representation and spectral density matrix of the multi-dimensional embedded self-similar process.\\

\begin{definition}
The process $U(t)=(U^0(t), U^1(t),\ldots,U^{q-1}(t))$ with parameter space
 $\check{T}=\{l^{n}, n\in {\Bbb Z}\}$ is a multi-dimensional self-similar process, where\\

$\bf (a)$\hspace{3mm} $\{U^j(\cdot)\}$ for every $j=0, 1, \cdots, q-1$ is self-similar process with parameter space

\hspace{12mm}$\check{T}=\{l^{n}, n\in {\Bbb Z}\}$.\\

$\bf (b)$\hspace{3mm} For every $n, \tau\in {\Bbb Z},\,\ j, k=0, 1, \cdots, q-1$
$$\mathrm{Cov}\big(U^j(l^{n+\tau}), U^k(l^{n})\big)=l^{2nH}\mathrm{Cov}\big(U^j(l^{\tau}),U^k(1)\big).$$
\end{definition}

\begin{remark}
Let $\{X(t), t \in\check{T}\}$ be a DSI process with scale $l>1$ and parameter space $\check{T}$, as defined in Remark {\em 2.1}.
Then by definition $3.1$ we have that  $\{U^j(l^n)\equiv X(l^n{\bf s}_j)\}$  for fixed $j=0, \ldots, q-1$, is a self-similar process and
$U(l^n)=\big(X(l^n{\bf s}_0), \ldots, X(l^n{\bf s}_{q-1})\big)$ where $1\leq {\bf s}_0<\cdots <{\bf s}_{q-1}<l$ is
a multi-dimensional self-similar process.
\end{remark}

By such method of sampling at discrete points, we
provide a $q$-dimensional embedded self-similar process $V(n)$ as
$$V(n)=\big(V^0(n), V^1(n), \ldots, V^{q-1}(n)\big),\hspace{1cm}n\in {\Bbb Z}$$
where $\{V^u(n)\equiv X(l^n{\bf s}_u)\}$ and $\{V(n)\equiv U(l^n)\}$
in which $U(l^n)$ and ${\bf s}_u$ follows the same assumptions as
in Remark 3.1.
Such definition of multi-dimensional embedded self-similar provides a platform
to obtain spectral density of $U(l^n)$ in the followings.

\begin{remark}
Corresponding to the $q$-dimensional embedded self-similar process
$V(n)$ there exists an embedded scale invariant process $\big
\{W(\kappa), \kappa \in {\Bbb Z}\big\}$ with scale $l$ as
\be \{W(\kappa)\equiv X(l^n{\bf s}_u)\}\hspace{1cm}\kappa\in {\Bbb Z}\ee
where $u=\kappa-q[\frac{\kappa}{q}]$, $n=[\frac{\kappa}{q}]$
and $\kappa=nq+u$, where by $(3.1)$
$$W(\kappa+q)\equiv X(l^{(n+1)}{\bf s}_u)\stackrel{d}{=}l^{H}X(l^n{\bf s}_u)\equiv l^{H}W(\kappa).$$
So the scale invariant process $\big \{ X(l^n{\bf s}_u), 1\leq
{\bf s}_0<\cdots <{\bf s}_{q-1}; n\in {\Bbb Z}\big\}$ and
the embedded scale invariant process $\big \{W(\kappa),\kappa\in{\Bbb Z}\big\}$, can be considered as counterparts.
\end{remark}

By the following theorem, the spectral representation and spectral density matrix of the $q$-dimensional embedded self-similar process and harmonic like representation of each column is obtained.

\begin{theorem}
Let $X(\cdot)$ be a {\em DSI} process with scale $l$ and $1\leqslant{\bf s}_0<\ldots<{\bf s}_{q-1}<l$,
then $V(n)=\big(V^0(n), \ldots, V^{q-1}(n)\big)$, where $\{V^u(n)\equiv X(l^n{\bf s}_u)\}$,
$n\in {\Bbb Z}$ and $u= 0, \ldots, q-1$ is a multi-dimensional embedded self-similar process and\\

\noindent {\em (i)} Harmonic like representation of $V^u(n)$ for fixed $u$
and $n\in {\Bbb Z}$ is
\be V^u(n)=(l^n{\bf s}_u)^H\int_{0}^{2\pi}e^{i\omega n}d\phi_{u}(\omega)\ee
where $\{\phi_u(\omega)\}$ are orthogonal spectral measures.
$E[d\phi_{u}(\omega)\overline{d\phi_{\nu}(\omega')}]=dG_{u,v}^H(\omega)$ when $\omega=\omega'$,
and is zero when $\omega\neq\omega'$ for $u, \nu= 0, \ldots, q-1$.\\

\noindent {\em (ii)} Spectral density matrix of $\{V(n)\equiv U(l^n)\}$ is
$g^H(\omega)=[g_{u,\nu}^H(\omega)]_{u,\nu=0, \ldots, q-1},$ where the elements $g_{u,v}^H(\omega)=dG_{u,v}^H(\omega)/d\omega$ are
\be g_{u,\nu}^H(\omega)=\frac{({\bf s}_u{\bf s}_{\nu})^{-H}}{2\pi}
\sum_{\tau=-\infty}^{\infty}l^{-H\tau}e^{-i\omega\tau}Q^H_{u,\nu}(\tau)\ee
$\tau\in{\Bbb N}$ and $Q^H_{u,\nu}(\tau)$ is the covariance function of $V^u(\tau)$ and $V^{\nu}(0)$.\\
\end{theorem}
{\bf Proof of (i):} Remark $2.4$ implies that
$$V^u(n)\equiv X(l^n{\bf s}_u)={\cal L}_{H}Y(l^n{\bf s}_u)=(l^n{\bf s}_u)^H\eta^u(n)$$
where $\eta^u(n)=Y(n\ln l+\ln{\bf s}_u)$. Thus $V^u(n)$ for every $u= 0, 1,\ldots, q-1$ is an embedded self-similar process in $n$,
where its discrete time stationary counterpart $\eta^u(n)$ for fixed $u= 0, 1,\ldots, q-1$ has spectral representation $\eta^u(n)=\int_{0}^{2\pi}e^{i\omega n}d\phi_{u}(\omega)$.\\ \\
{\bf Proof of (ii):} The covariance matrix of $V(n)$ is denoted by $Q^H(n,\tau)=[Q^H_{u,\nu}(n, \tau)]_{u,\nu=0, \ldots, q-1}$ where
$$Q^H_{u,\nu}(n, \tau)=E[V^u(n+\tau)V^{\nu}(n)]=E[X(l^{(n+\tau)}{\bf s}_u)X(l^n{\bf s}_{\nu})]$$\\
By the scale invariant property of the process $X(\cdot)$ we have that\\
\be Q^H_{u,\nu}(n, \tau)=l^{2nH}E[X(l^{\tau}{\bf s}_u)X({\bf s}_{\nu})]=l^{2nH}Q^H_{u,\nu}(\tau)\ee\\
where $Q^H_{u,\nu}(\tau)=Q^H_{u,\nu}(0, \tau)=E[V^u(\tau)\overline{V^{\nu}(0)}]$, then by (3.2)
$$Q^H_{u,\nu}(\tau)=E[(l^{\tau}{\bf s}_u)^H({\bf s}_{\nu})^H\int_{0}^{2\pi}e^{i\omega\tau}d\phi_{u}(\omega)
\int_{0}^{2\pi}\overline{d\phi_{v}(\omega')}]$$
\be =l^{\tau H}({\bf s}_u{\bf s}_{\nu})^H\int_{0}^{2\pi}e^{i\omega\tau}dG^H_{u,\nu}(\omega)\ee\\
where $E[d\phi_{u}(\omega)\overline{d\phi_{\nu}(\omega')}]=dG^H_{u,\nu}(\omega)$ when $\omega=\omega'$ and is $0$ when $\omega\neq\omega'$.\\
On the other hand, by the definition of $\eta^u(n)$ in the proof of part $(i)$\\
$$Q^H_{u,\nu}(\tau)=E[X(l^{\tau}{\bf s}_u)X({\bf s}_{\nu})]
=E[{\cal L}_{H}Y(l^{\tau}{\bf s}_u){\cal L}_{H}Y({\bf s}_{\nu})]$$
$$=(l^{\tau}{\bf s}_u{\bf s}_{\nu})^HE[Y(\tau\ln l+\ln{\bf s}_u)Y(\ln{\bf s}_{\nu})]$$
$$=(l^{\tau}{\bf s}_u{\bf s}_{\nu})^HE[\eta^u(\tau)\eta^{\nu}(0)]
=(l^{\tau}{\bf s}_u{\bf s}_{\nu})^HB_{u, \nu}(\tau).$$
As $\eta^u(\cdot)$ ia a stationary process so by (2.3)
$$B_{u, \nu}(\tau)=\int_{0}^{2\pi}e^{i\omega\tau}dG^H_{u,\nu}(\omega),\hspace{5mm}u, \nu= 0, \ldots, q-1$$
Now by (2.4) for $u,\nu=0, \ldots, q-1$ we have
$$G^H_{u,\nu}(A)=\frac{1}{2\pi}\int_{A}\sum_{\tau=-\infty}^{\infty}B_{u, \nu}(\tau)e^{-i\lambda\tau}d\lambda.$$
By substituting $B_{u, \nu}(\tau)=(l^{\tau}{\bf s}_u{\bf s}_{\nu})^{-H}Q^H_{u,\nu}(\tau)$, the elements of the spectral distribution function, $G^H_{u,\nu}(\cdot)$ has the following representation\\
\be G^H_{u,\nu}(A)=\frac{({\bf s}_u{\bf s}_{\nu})^{-H}}{2\pi}
\int_{A}\sum_{\tau=-\infty}^{\infty}l^{-H\tau}e^{-i\lambda\tau}Q^H_{u,\nu}(\tau)d\lambda.\ee

Let $A=(\omega,\omega+d\omega]$, then the elements of the spectral density matrix, $g_{u,\nu}^H(\omega)$ are
$$g_{u,\nu}^H(\omega):=\frac{G^H_{u,\nu}(d\omega)}{d\omega}=\frac{({\bf s}_u{\bf s}_{\nu})^{-H}}{2\pi}\sum_{\tau=-\infty}^{\infty}l^{-H\tau}\big(\frac{1}{d\omega}
\int_{\omega}^{\omega+d\omega}e^{-i\lambda\tau}d\lambda\big)Q^H_{u,\nu}(\tau)$$
$$=\frac{({\bf s}_u{\bf s}_{\nu})^{-H}}{2\pi}\sum_{\tau=-\infty}^{\infty}l^{-H\tau}
\big(\frac{1}{-i\tau}\lim_{d\omega\rightarrow 0}
\frac{e^{-i{(\omega+d\omega)\tau}}-e^{-i\omega \tau}}{d\omega}\big)Q^H_{u,\nu}(\tau).$$

Thus we get to the assertion of part (ii) of the theorem.$\square$\\

\renewcommand{\theequation}{\arabic{section}.\arabic{equation}}
\section{Scale invariant process with Markov property}
\setcounter{equation}{0}
Using the method of sampling in section 3 and the scale invariant process $X(\cdot)$ in Remark 3.1, and
its corresponding embedded scale invariant process $W(\cdot)$ in Remark 3.2, we introduce $W(\cdot)$ with Markov property in the wide sense, named embedded scale invariant Markov process. Also the  corresponding multi-dimensional embedded self-similar process $V(\cdot)$ and its corresponding multi-dimensional self-similar process $U(\cdot)$ is defined. We find the covariance function of embedded scale invariant Markov process in subsection 4.1. The spectral density matrix of multi-dimensional process of it embedded is evaluated in subsection 4.2.

\subsection{Covariance function of embedded scale invariant Markov process}
Here we characterize the covariance function of the embedded scale invariant Markov process $\{W(\kappa), \kappa\in {\Bbb Z}\}$ in Theorem 4.1 and the covariance function of the associated multi-dimensional embedded self-similar Markov process in Theorem 4.2.
\begin{theorem}
Let $\{W(\kappa), \kappa\in{\Bbb Z}\}$, defined by $(3.1)$, be an embedded scale invariant and Markov process in the wide sense with scale $l$. Then for $\tau\in{\Bbb W}$, $\kappa=nq+\nu$, $\kappa+\tau=mq+u$, $u, \nu= 0, \ldots, q-1$ and $n,m \in {\Bbb Z}$, the covariance function
\be R_{\kappa}(\tau):=E[W(\kappa+\tau)W(\kappa)]=E[X(l^{m}{\bf s}_u)X(l^n{\bf s}_{\nu})]\ee
where $1\leqslant{\bf s}_0< {\bf s}_1<\ldots<{\bf s}_{q-1}<l$, can be characterized as
\be R_{\kappa}(tq+w)=[\tilde{f}(q-1)]^t\tilde{f}(\kappa+w-1)[\tilde{f}(\kappa-1)]^{-1}R_{\kappa}(0),\ee
$t\in{\Bbb Z}$, $w=0, \ldots, q-1$
\be\tilde{f}(r)=\prod_{j=0}^{r}f(j)=\prod_{j=0}^{r}R_j(1)/R_j(0),\hspace{7mm}r\in{\Bbb Z}, \tilde{f}(-1)=1.\ee
\end{theorem}

\noindent
{\bf Proof:} By considering sample points of this paper,  proof of this theorem follows by a similar method as for theorem 3.2 of [6].

Now we can use this theorem to prove the next result for multi-dimensional embedded self-similar Markov process.

\begin{theorem}
Let $\{W(\kappa),\kappa\in {\Bbb Z}\}$ be a embedded scale invariant Markov process,
and $\{V(n)\equiv U(l^n)\}$ be its associated multi-dimensional embedded self-similar Markov, where $\{U(l^n), n\in{\Bbb Z}\}$ is the corresponding self-similar Markov process, both with  the same covariance matrix $Q^H(n, \tau)$ which is defined by
$(3.5)$. Then\\
\be Q^H(n,\tau)=l^{2nH}[\tilde{f}(q-1)]^{\tau}CD,\hspace{7mm}\tau\in{\Bbb Z}\ee
where $\tilde{f}(\cdot)$ is defined in $(4.3)$ and the matrices $C$ and $D$ are given by $C=[C_{u, \nu}]_{u,\nu= 0, \ldots, q-1}$, where $C_{u,\nu}=\tilde{f}(u-1)[\tilde{f}(\nu-1)]^{-1}$, and $D$ is a diagonal matrix with diagonal elements $R_{\nu}(0)$, $\nu=0, 1, \ldots, q-1$, which is defined in $(4.1)$.\\
\end{theorem}
{\bf Proof:} As $W(\cdot)$ is embedded scale invariant with scale $l$, (3.2) and (3.5) indicate that $Q^H_{u,\nu}(n, \tau)=l^{2nH}Q^H_{u,\nu}(\tau)$. Now by the assumption $\kappa=nq+\nu$ and $\kappa+\tau=mq+u$ where $m,n\in {\Bbb Z}$, $\tau\in{\Bbb W}$, we have $\tau=(m-n)q+u-\nu$ and therefore
$$R_{\kappa}(\tau)=R_{nq+\nu}((m-n)q+u-\nu)=E[W(mq+u)W(nq+\nu)]$$
$$=E[X(l^{m}{\bf s}_u)X(l^n{\bf s}_{\nu})].$$
Hence
\be Q^H_{u,\nu}(\tau)=E[X(l^{\tau}{\bf s}_u)X({\bf s}_{\nu})]=R_{\nu}(\tau q+u-\nu)\ee\\
and by the Markov property of $W(\cdot)$ we have
$$R_{\nu}(\tau q+u-\nu)=[\tilde{f}(q-1)]^{\tau}\tilde{f}(u-1)[\tilde{f}(\nu-1)]^{-1}R_{\nu}(0)$$
for $u, \nu= 0, \ldots, q-1$. Let $C_{u, \nu}=\tilde{f}(u-1)[\tilde{f}(\nu-1)]^{-1}$, so\\
\be Q^H_{u,\nu}(\tau)=[\tilde{f}(q-1)]^{\tau}C_{u, \nu}R_{\nu}(0).\ee
Thus we can represent the elements of the covariance matrix of $q$-dimensional embedded self-similar Markov process as
$$Q^H_{u,\nu}(n, \tau)=l^{2nH}[\tilde{f}(q-1)]^{\tau}C_{u, \nu}R_{\nu}(0).\square$$

\subsection{Spectral representation of the process}
The spectral density matrix of the multi-dimensional embedded self-similar Markov
and corresponding multi-dimensional self-similar Markov processes are characterized by the
following proposition.

\begin{proposition}
The spectral density matrix $g^H(\omega)=[g^H_{u,\nu}(\omega)]_{u,\nu= 0, \ldots, q-1}$ of the $q$-dimensional embedded self-similar Markov process $\{V(n)\equiv U(l^n)\}$, where $U(l^n)$ is the corresponding $q$-dimensional self-similar Markov
process, is specified by
$$g_{u,\nu}^H(\omega)=\frac{({\bf s}_u{\bf s}_{\nu})^{-H}}{2\pi}
\left[\frac{\tilde{f}(u-1)R_{\nu}(0)}{\tilde{f}(\nu-1)
\big(1-e^{-i\omega}l^{-H}\tilde{f}(q-1)\big)}-
\frac{\tilde{f}(\nu-1)R_u(0)}{\tilde{f}(u-1)\big(1-e^{-i\omega}l^{H}
\tilde{f}^{-1}(q-1)\big )}\right]$$
where $R_{k}(0)$ is the variance of $W(k)$ and $\tilde{f}(\cdot)$ is defined by $(4.3)$.
\end{proposition}
{\bf Proof:} By applying (3.3) and (4.12), the spectral density matrix of the
 process $\{V(n), n\in {\Bbb Z}\}$ which is denoted by $g^H(\omega)=[g^H_{u, \nu}(\omega)]_{u,\nu= 0, \ldots, q-1}$ can be written as
$$g_{u,\nu}^H(\omega)=\frac{({\bf s}_u{\bf s}_{\nu})^{-H}}{2\pi}
\Big[\sum_{\tau=0}^{\infty}l^{-H\tau}e^{-i\omega\tau}Q^H_{u,\nu}(\tau)$$
\vspace{-5mm}
$$+\sum_{\tau=-\infty}^{-1}l^{-H\tau}e^{-i\omega\tau}Q^H_{u,\nu}(\tau)\Big]=g_{u,\nu,1}^H(\omega)+g_{u,\nu,2}^H(\omega)$$
where
$$g_{u,\nu,1}^H(\omega)=\frac{({\bf s}_u{\bf s}_{\nu})^{-H}}{2\pi}
\sum_{\tau=0}^{\infty}l^{-H\tau}e^{-i\omega\tau}[\tilde{f}(q-1)]^{\tau}\tilde{f}(u-1)[\tilde{f}(\nu-1)]^{-1}R_{\nu}(0)$$
\be =\frac{({\bf s}_u{\bf s}_{\nu})^{-H}\tilde{f}(u-1)R_{\nu}(0)}{2\pi\tilde{f}(\nu-1)}
\sum_{\tau=0}^{\infty}\big(l^{-H}e^{-i\omega}\tilde{f}(q-1)\big)^{\tau}.\ee\\
By Remark 3.2, the scale invariant property of $W(\kappa)$ and the assumption, that at least one of the $\text{Corr}[W(j)W(j+1)]$ be smaller than one, we have that $|\tilde{f}(q-1)|<l^{H}$ for $j= 0, \ldots, q-1$. Thus
$$|e^{-i\omega}l^{-H}\tilde{f}(q-1)|=|l^{-H}\tilde{f}(q-1)|<1,$$
and (4.13) for $\tau\in{\Bbb W}$ is convergent. By the equality
$$Q_{u,\nu}(-\tau)=E[X(l^{-\tau}{\bf s}_u)X({\bf s}_{\nu})]=l^{-2\tau H}E[X(l^{\tau}{\bf s}_{\nu})X({\bf s}_{u})]
=l^{-2\tau H}Q_{\nu, u}(\tau),$$
convergence of $g_{u,\nu,2}^H(\omega)$ follows by a similar method. Therefore\\
$$g_{u,\nu}^H(\omega)=\frac{({\bf s}_u{\bf s}_{\nu})^{-H}}{2\pi}
\Big[\frac{R_{\nu}(0)\tilde{f}(u-1)}{\tilde{f}(\nu-1)}\sum_{\tau=0}^{\infty}\big(l^{-H}e^{-i\omega}\tilde{f}(q-1)\big)^{\tau}$$
$$+\frac{R_{u}(0)\tilde{f}(\nu-1)}{\tilde{f}(u-1)}\sum_{\tau=1}^{\infty}\big(l^{-H}e^{i\omega}\tilde{f}(q-1)\big)^{\tau}\Big]$$
$$=\frac{({\bf s}_u{\bf s}_{\nu})^{-H}}{2\pi}\Big[\frac{R_{\nu}(0)\tilde{f}(u-1)}{\tilde{f}(\nu-1)\big(1-l^{-H}e^{-i\omega}\tilde{f}(q-1)\big)}
+\frac{R_{u}(0)\tilde{f}(\nu-1)l^{-H}e^{i\omega}\tilde{f}(q-1)}{\tilde{f}(u-1)\big(1-l^{-H}e^{i\omega}
\tilde{f}(q-1)\big)}\Big],$$\\
so we arrive at the conclusion of the proposition.$\square$\\

\begin{example}
Let
\be X(t)=\sum_{n=1}^{\infty}\lambda^{n(H-\frac{1}{2})}I_{[\lambda^{n-1}, \lambda^n)}(t)B(t)\ee
where $B(\cdot)$  is the standard  Brownian motion, $I(\cdot)$ indicator function, $H>0$ and  $\lambda>1$. This process is a Brownian motion inside each scale $[\lambda^{n-1}, \lambda^n)$ and in general is a {\em DSI} process with scale $\lambda$ and Hurst index $H$.
For $H=0.5$ this process is just standard Brownian motion, which is a scale invariant process with Hurst index $H$. For $H\neq 0.5$, we call $X(t)$ a Simple Brownian Motion ({\em SBM}).
We showed in {\em \cite{m2}} that $\{X(t), t\in {\Bbb R}^+\}$ is {\em DSI} and Markov  with Hurst index $H$ and scale $ \lambda $.
By sampling of this process  at points $\lambda^n{\bf s}_u$, $n\in {\Bbb W}$, where
$1\leq {\bf s}_0<\ldots <{\bf s}_{q-1}<\lambda$, and  by assuming $\lambda=\alpha^T$, we have the  corresponding multi-dimensional self-similar process as $U(\lambda^{n})=\big(X(\lambda^{n}{\bf s}_0), \ldots, X(\lambda^{n}{\bf s}_{q-1})\big)$.
So
$\{W(\kappa )\equiv X(\lambda^{n}{\bf s}_u)\},$ is an embedded scale invariant Markov process,  and  $V(n)=\big(V^0(n), \ldots, V^{q-1}(n)\big)$ where $\{V^u(n)\equiv W(\kappa)\}$ is the associated $q$-dimensional embedded self-similar Markov process
where $u=\kappa-q[\frac{\kappa}{q}]$, $n=[\frac{\kappa}{q}]$.
By {\em (4.1)} we have that  $R_j^H(0)=R_ j^H(1)=\lambda^{2H'}{\bf s}_j$
for $j=0,\cdots ,q-2$ and $R_{q-1}^H(1)=\lambda^{H'}R_{q-1}^H(0)=\lambda^{3H'}{\bf s}_{q-1} $, where $H'= H-\frac{1}{2}$. So
$R_u(0)=\lambda^{2H'}{\bf s}_u$,
$\;R_{\nu}(0)=\lambda^{2H'}{\bf s}_{\nu}$. Also {\em (4.3)} implies  that
$\tilde{f}(u-1)=\tilde{f}(\nu-1)=1$, $\tilde{f}(q-1)=\lambda^{H'}$. Thus  By Proposition {\em 4.1}, the spectral density matrix of $V(n)$ is given by $g^H(\omega)$ where
$$g_{u,\nu}^H(\omega)=\frac{({\bf s}_u{\bf s}_{\nu})^{-H}\lambda^{2H'}}{2\pi}\left[
\frac{{\bf s}_{\nu}}{1-e^{-i\omega}\lambda^{-1/2}}-\frac{{\bf s}_{u}}{1-e^{-i\omega}\lambda^{1/2}} \right].$$
\end{example}

\section{Simulation}
{We have used  Matlab program to  simulate and plot SBM defined by (4.14) and its corresponding multi-dimensional
self-similar process for different values of $H$ and $\lambda $. We have simulated
$$X(t)=\sum_{n=1}^{M}\lambda^{n(H-\frac{1}{2})}I_{[\lambda^{n-1},\lambda^n)}(t)B(t)$$
where $M=30$. We also assume to have $q=30$ samples  in each scale interval $I_n:=[\lambda^{n-1}, \lambda^n)$,
where $n=1,2,\cdots, M$. By choosing these sample points to be  $\lambda^{n-1}{\bf s}_i$, where ${\bf s}_i=1+i(a-1)/T$
for $i=0,1, \cdots, q-1$, our sample points will be equally spaced in each scale interval $I_n,\; n=1,2,\cdots, M$.\\

\input{epsf}
\epsfxsize=7.8cm \epsfysize=4cm
\begin{figure}
\hspace{-.7cm}{\centering\epsffile{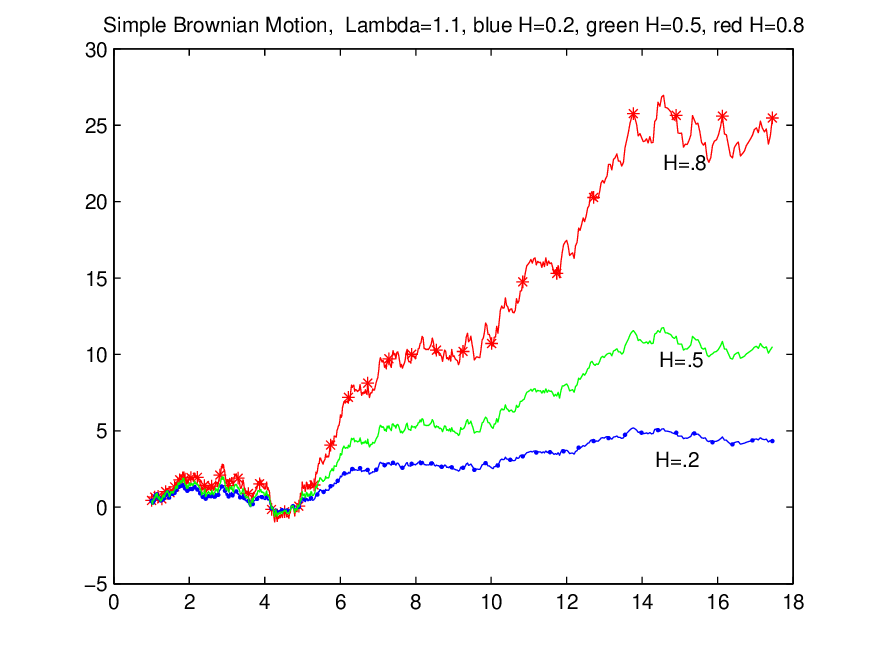} \hspace{-.9cm}
 \epsfxsize=7.8cm \epsfysize=4cm \epsffile{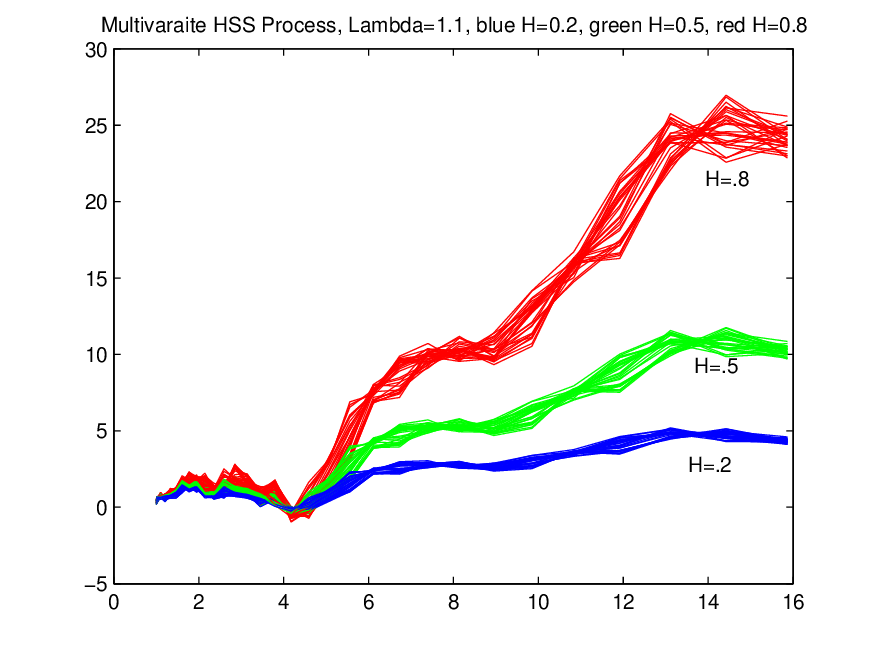 }\epsfxsize=4cm \epsfysize=4cm \vspace{-.3cm}
\caption{\tiny Simple Brownian Motion and corresponding Multivariate
H-ss Process where Scales and Hurst index are indicated
.$\hspace{1.7cm}$above the figures  }}$\vspace{-0.2in}$\\
\end{figure}

\input{epsf}
\epsfxsize=7.8cm \epsfysize=4cm
\begin{figure}
\hspace{-.7cm}{\centering\epsffile{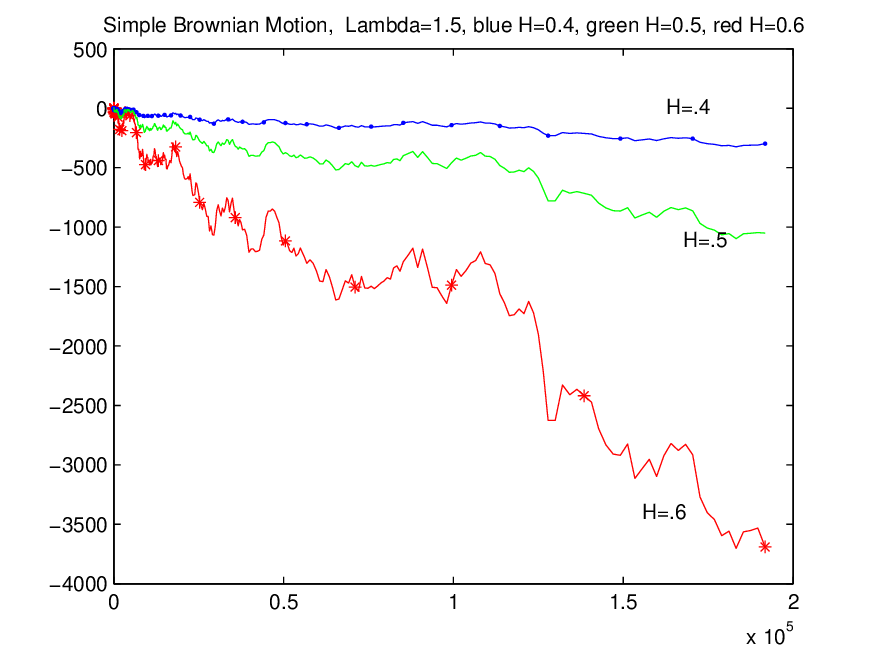} \hspace{-.9cm}
 \epsfxsize=7.8cm \epsfysize=4cm \epsffile{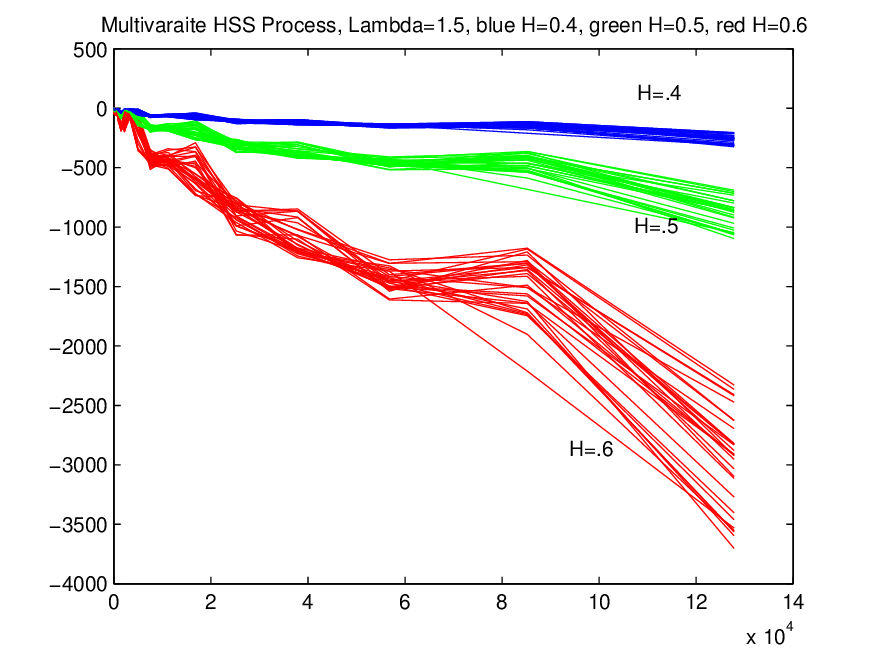 }\epsfxsize=4cm \epsfysize=4cm\vspace{-.3cm}
 \caption{\tiny Simple Brownian Motion and
corresponding Multivariate H-ss Process where Scales and Hurst index
are indicated .$\hspace{1.7cm}$ above the figures}}
\end{figure}
All the multi-dimensional self-similar processes have been plotted at points $\lambda^n$ where $n=0,1,2, \cdots, M$.
One can easily verify that how SBM's  are going to enlarge  at the beginning of each scale interval $[\lambda^{n-1}, \lambda^n)$, which is the main property of DSI processes, while for the corresponding
multi-dimensional self-similar process, which have been constructed
by  one observation in each scale interval, these jumps are equally
like  for all observations. Figure 1, consists of two figures SBM on
the left and corresponding multi-dimensional self-similar on the
right. The figure on the left consists of three different curves of
SBM, all with scale $\lambda=1.1$, but with different Hurst index.
It is worthy to note that we have simulated just one discrete time
Brownian motion $B(\lambda^n{\bf s}_i)$ of Example 4.1, for these three
curves. The curve in the middle  has Hurst index $H=0.5$, so it is a
discrete time Brownian motion which is a self-similar process and
other two curves are to compare with this. The upper curve has Hurst
index $H=0.8$, so it is a scale invariant process and at the beginning of each
scale interval $[\lambda^{n-1}, \lambda^n)$  enlargement in compare
with  Brownian motion occurs, which has been caused by the growth of
coefficients to $ \lambda^{n(0.8-1/2)}$  at the beginning of $n$-th
scale interval. Also the lower  curve has scale $H=0.2$, so in
compare with Brownian motion,  the coefficient at the beginning of $n$-th scale interval decreases
to $\lambda^{n(0.2-1/2)}$. So it comes to have less variation than Brownian Motion
at the beginning of each scale interval. Figure 2 is also included
of two figures, where the left one again  consist of three curves of SBM  all with scale
$\lambda=1.5$, but with Hurst indices $H=0.4, 0.5$ and $0.6$. Again we have generated one $B(t)$ for all these three curves.
The curve in the middle has Hurst index $H=0.5$, so is a Brownian motion.
The upper curve is SBM with $H=0.6$, where enlargement in compare to Brownian motion occurs at the beginning of scale intervals
by $\lambda^{n(0.6-1/2)}$ and  to the same direction of the Brownian motion, the middle curve.
The curve with   $H=0.4$, lower curve , has  less variation in compare with  Brownian motion. So the size of lines  at
the beginning of all scale intervals decreases,by the rate
$\lambda^{n(0.4-1/2)}$ with respect to the Brownian motion, for the $n$-th scale interval.
For the corresponding  multi-dimensional self-similar processes where the $i$-th curve, for $i=0,\ldots, q-1$
has been evaluated at sample points $\big\{\lambda^n{\bf s}_i, n\in {\Bbb Z}\big\}$
of corresponding SBM, and has been plotted at points $\lambda^n$, and all are self-similar with the same Hurst index.
Finally the growth of Hurst index causes the growth of all lines in the corresponding multi-dimensional self-similar process as well.

One can compare these multi-dimensional self-similar  processes with the corresponding SBM which shows that as these are close together at any point $\lambda^n$, the changes inside the corresponding scale interval are less. It is also interesting that as at the end points of the multivariate self-similar process for $H=0.8$, have more variation , so the variation inside last scale intervals of SBM with $H=0.8$ is more. Finally as the path of all multi-dimensional self-similar for $H=0.8$ and $\lambda=1.5$ for last observations are decreasing, and for $H=0.8$ and $\lambda=1.1$ are increasing, so the path of corresponding SBM at the beginning of last scale intervals are respectively decreasing and increasing.
\hspace{-5in} \nopagebreak


\end{document}